\documentclass{aptpub}
\authornames{MAMADOU BA, AHMADOU BAMBA SOW AND ETIENNE PARDOUX} % insert the authors here for use in running head
\shorttitle{Binary trees, exploration processes, and an extented Ray--Knight Theorem} % insert short title here for use in running head

\numberwithin{equation}{section}

\newcommand{\eps}  {\varepsilon}

\usepackage{rotating}
\usepackage{amsmath, amssymb, amsfonts}
\newcommand\unnumberedfootnote[1]{ %
        \let\temp=\thefootnote %
        \renewcommand{\thefootnote}{}%
        \footnote{#1}%
        \let\thefootnote=\temp%
        \addtocounter{footnote}{-1}}

\newcommand{\R}{\mathbb{R}}
\renewcommand{\L}{\mathcal{L}}
\newcommand{\sign}  {\text{sign}}

\def \bpf {\noindent{\sc Proof: }}
\def \epf {\hbox{}\nobreak\hfill
\vrule width 2mm height 2mm depth 0mm
\par \goodbreak \smallskip}

\usepackage{graphicx}
\usepackage{graphics}
\usepackage{pstricks}
\usepackage{subfigure}
\usepackage{color}
\definecolor{DarkGreen}{rgb}{0,0.6,0}
\definecolor{VertFonce}{rgb}{0,0.6,0}
\definecolor{MidGreen}{rgb}{0.6,1,0.6}
\definecolor{LightGreen}{rgb}{0.88,1,0.88}
\definecolor{LightGray}{rgb}{0.94,0.94,0.94}
\definecolor{VeryLightBlue}{rgb}{0.9,0.9,1}
\definecolor{LightBlue}{rgb}{0.8,0.8,1}
\definecolor{DarkBlue}{rgb}{0,0,0.6}
\definecolor{VeryLightYellow}{rgb}{1,1,0.9}
\definecolor{LightYellow}{rgb}{1,1,0.6}
\definecolor{MidYellow}{rgb}{1,1,0.5}
\definecolor{VeryLightRed}{rgb}{1,0.9,0.9}
\definecolor{LightRed}{rgb}{1,0.8,0.8}
\definecolor{turquoise}{rgb}{0.00,0.53,0.68}{}
\definecolor{mauve}{rgb}{0.50,0.00,0.50}

\begin{document}
\title{Binary trees, exploration processes, and an extented Ray--Knight Theorem}

\authorone[LATP--Universit\'e de Provence]{Mamadou Ba}
\addressone{Centre de Math\'{e}matiques et d'Informatique, 
39 rue F. Joliot-Curie, 13453, Marseille, France. email: mba@cmi.univ-mrs.fr} 
\authortwo[LATP--Universit\'e de Provence]{Etienne Pardoux}
\addresstwo{Centre de Math\'{e}matiques et d'Informatique, 
39 rue F. Joliot-Curie, 13453, Marseille, France. email: pardoux@cmi.univ-mrs.fr} 
\authorthree[LERSTAD, UFR S.A.T., U.G.B.]{Ahmadou Bamba Sow}
\addressthree{Universit\'e Gaston Berger, BP 234, Saint-Louis, Senegal. email :
ahmadou-bamba.sow@ugb.edu.sn}

\begin{abstract}
{\small 
We study the bijection between binary Galton--Watson trees in continuous time and their exploration process, both in the sub- and in the supercritical cases. We then take the limit over renormalized quantities, as the size of the population tends to infinity. We thus deduce Delmas' generalization of the second Ray--Knight theorem.
 }
\end{abstract}
\keywords{Galton--Watson processes, Feller's branching, Exploration process, Ray--Knight theorem.}
\ams{60J80;60F17} {92D25}

 \section*{Introduction}
 There are various forms of bijection between an exploration (or height) process and a random binary tree. This paper starts with a description of such a bijection, and a new rather simple proof that a certain law on the exploration paths is in bijection with the law of a continuous time binary Galton--Watson random tree. The result in the critical case was first established by Le Gall \cite{JFLG}, and in the subcritical case by Pitman and Winkel \cite{PiWi}, see also Geiger and Kersting \cite{GeKe}, Lambert \cite{Lam}, where the exploration processes are jump processes, while ours are continuous.  For similar results in the case where the approximating process is in discrete time and the tree is not necessarily binary, see Duquesne and Le Gall \cite{DLG}.
 We consider also the supercritical case, which is new. Inspired by the work of Delmas \cite{JFD}, we note that in the supercritical case, the random tree killed at time $a>0$ is in bijection with the exploration process reflected below $a$. Moreover, one can define a unique local time process, which describes the local times of all the reflected exploration processes, and has the same law as the supercritical Galton--Watson process. 
  We next renormalize our Galton--Watson tree and height process, and take the weak limit, thus providing a new proof of Delmas' extension \cite{JFD} of the second Ray--Knight theorem. The classical version of this theorem establishes the identity in law 
 between the local time of reflected Brownian motion considered at the time when the local time at $0$ reaches $x$, and at all levels, and a Feller critical branching diffusion. The same result holds in the subcritical (resp. supercritical) case, Brownian motion being replaced by Brownian motion with drift (in the supercritical case, reflection below an arbitrary level,
 as above, is needed).  The exploration process in fact describes the genealogical tree (in the sense of Aldous \cite{ALD}) of the population, whose total mass follows a Feller SDE. Our proof by approximation makes this interpretation completely transparent. 
 
 The paper is organized as follows. Section 1 is devoted to the description of the bijection between height curves and binary trees.  Section 2 presents the relation between laws of height processes and Galton--Watson trees, and the ``discrete Ray Knight theorem''.  Section 3 presents the results of convergence of both the population process and the height process, in the limit of large populations. Finally section 4 deduces the generalized Ray--Knight theorem from our convergences and the results at the discrete level.
 
 \section{Preliminaries}\label{prelim}
 Fix $p>0$. 
 Consider a continuous piecewise linear function $H$ from a subinterval of $\R_+$ into $\R_+$, which possesses the following properties~: 
its slope is either $p$ or $-p$; it starts at time $t=0$ from 0 with the slope $p$; whenever $H(t)=0$, $H'_-(t)=-p$ and
$H'_+(t)=p$; $H$ is
 stopped at the time $T_m$ of its $m$--th return to $0$, which is supposed to be finite. We shall denote by $\mathcal{H}_{p,m}$
the collection of all such functions.
 We shall write $\mathcal{H}_p$ for $\mathcal{H}_{p,1}$. We add the restriction that between two consecutive visits to zero, any function from $\mathcal{H}_{p,m}$ has all its local minima at distinct heights.
 
We denote by $\mathcal{T}$ the set of finite rooted binary trees which are defined as follows. An ancestor is born at time 0. Until she eventually dies, she gives birth to an arbitrary number of offsprings, but only one at a time. 
 The same happens to each of her offsprings, and the offsprings of her offsprings, etc... until eventually the population dies out. 
 We denote by $\mathcal{T}_m$ the set of forests which are the union of $m$ elements of $\mathcal{T}$.
 
\begin{figure}[htbp]
 \centering
 \includegraphics[scale=0.6]{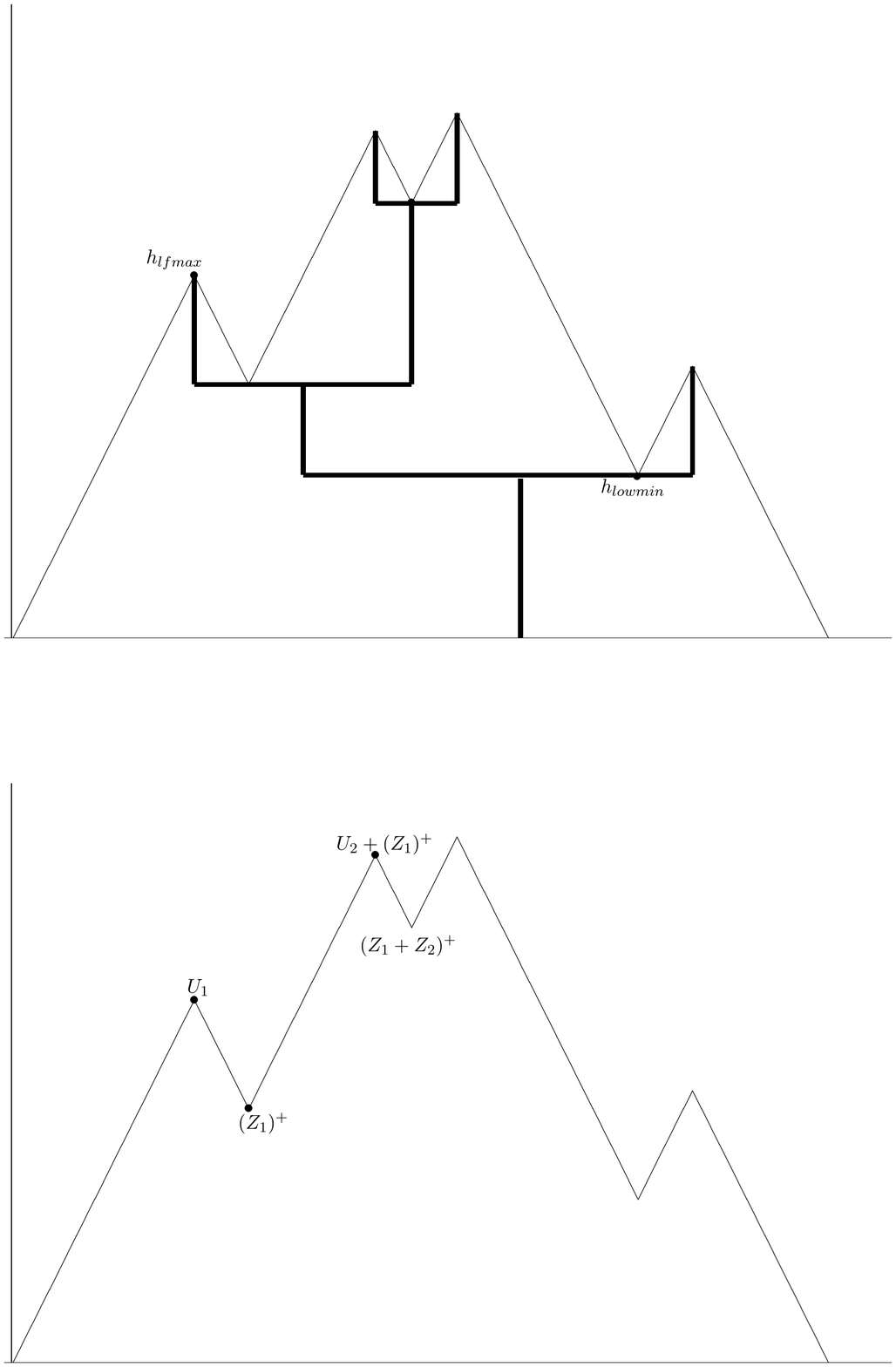} 
 \caption {\label{Fig 1} Bijection between $\mathcal{H}_2$ and $\mathcal{T}$ (see above), and a trajectory of an exploration process (see below)}
 \end{figure}

There is a well known bijection between binary trees and exploration processes. Under the curve representing an element $H\in \mathcal{H}_p$, we can draw a tree as follows. The height $h_{lfmax}$ of the leftmost local maximum of $H$ is the lifetime of the ancestor and the height $h_{lowmin}$ of the lowest non zero local minimum is the time of the birth of the first offspring of the ancestor. If there is no such local minimum, the ancestor dies before giving birth to any offspring. We draw a horizontal line at level $h_{lowmin}$. $H$ has two excursions above $h_{lowmin}$. The right one is used to represent the fate of the first offspring and of her progeny. The left one is used to represent the fate of the ancestor and of the rest of her progeny, excluding the first offspring and her progeny. If there is no other local minimum of $H$ neither at the left, nor at the right of the first explored one, it means that there is no further birth: we draw a vertical line up to the unique local maximum, whose height is a death time. Continuing until there is no further local minimum-maximum to explore, we
 define by this procedure a bijection $\Phi_p$ from $\mathcal{H}_p$ into $\mathcal{T}$ (see Figure \ref{Fig 1}). 
 Repeating the same construction $m$ times, we extend $\Phi_p$ as a bijection between $\mathcal{H}_{p,m}$ and $\mathcal{T}_m$. Note that describing the exploration process from a tree is obvious (the horizontal distances between the vertical branches are arbitrary). See the top of Figure 1. 
 
We now define probability measures on $\mathcal{H}_p$ (resp. $\mathcal{H}_{p,m}$) and $\mathcal{T}$ (resp. $\mathcal{T}_m$). 
We describe first the subcritical case (by a slight abuse of terminology, subcritical in the present paper always  means either subcritical  or critical). Let $0<\mu\le\lambda$ be two parameters. We define a stochastic process whose trajectories belong to $\mathcal{H}_p$ as follows.
 Let $\left\lbrace U_k ,~k\geq 1\right\rbrace$ and $\left\lbrace
V_k,~k\geq 1 \right\rbrace$ be two mutually independent sequences of i.i.d exponential random variables with means  $1/\lambda$ and  $1/\mu$ respectively. Let $Z_k=U_k-V_k$, $k\ge1$. $\mathbb{P}_{\lambda,\mu}$ is the law of the random element of $\mathcal{H}_p$, which is such that the height of the first local  maximum is $U_1$, that of the first local minimum is $(Z_1)^+$.  If $(Z_1)^+=0$, the process is stopped. Otherwise, the height of the second local maximum is $Z_1+U_2$, the height of the second local minimum is $(Z_1+Z_2)^+$, etc. Because $\mu\le\lambda$, the process returns to zero a.s. in finite time. The random trajectory which we have constructed is  an excursion above zero (see the bottom of Figure \ref{Fig 1}). We define similarly a law on $\mathcal{H}_{p,m}$ as the concatenation of $m$ i. i. d. such excursions, and denote it by $\mathbb{P}_{\lambda,\mu}$. This random element defined above is called an exploration process or height process.
We associate the continuous time Galton--Watson  tree (which is a random element of $\mathcal{T}$) with the same pair of parameters $(\lambda,\mu)$ as follows. The lifetime of each individual is exponential with expectation $1/\lambda$, and during her lifetime, independently of it, each individual gives birth to offsprings according to a  Poisson process with rate $\mu$.  The behaviors of the various individuals are i. i. d.  We denote by $\mathbb{Q}_{\lambda,\mu}$ the law on  $\mathcal{T}_m$ of 
a forest of $m$ i. i d. random trees whose law is as just described. 

In the supercritical case, the case where $\lambda > \mu$, the exploration process defined above does not come back to zero $a.s.$. To overcome this difficulty, we add reflection below an arbitrary level $a>0$, and  we consider the height process $H^a=\{H^a_t,\ t\ge0\}$ reflected at level $a$ defined as above, with the addition of the rule that whenever the process reaches the level $a$, it stops and starts immediately going down with slope $-p$ for a duration of time exponential with expectation $1/\mu$. Again the process stops when first going back to zero. The reflected process $H^a$ comes back to zero almost surely. Indeed, let $A^a_n$ denote the event  ``$H^a$ does not reach zero during its $n$ first descents''.
We have clearly, since the levels of local maxima are bounded by $a$,
$\mathbb{P}(A_n)\leq (1-\exp(-\mu a))^n$, which goes to  zero as $n\longrightarrow \infty $. Hence the result. For each $a\in (0,+\infty)$, and any pair $(\lambda,\mu)$ of positive numbers, denote by $\mathbb{P}_{\lambda,\mu,a}$ the law of the process $H^a$. Define $\mathbb{Q}_{\lambda,\mu,a}$ to be the law of a binary Galton--Watson tree with birth rate $\mu$ and death rate $\lambda$, killed at time $t=a$ (i. e. all individuals alive at time $a^-$ are killed at time $a$). $\mathbb{P}_{\lambda,\mu,+\infty}$ makes perfect sense in case $\mu\leq\lambda$, $\mathbb{Q}_{\lambda,\mu,+\infty}$ is always well defined.

\section{Correspondence of laws}
The aim of this section is to prove that, for any  $\lambda,\mu >0$ and $a\in (0,+\infty)$ [including the case $a=+\infty$ in the case $\mu\leq\lambda$] ,  $\mathbb{P}_{\lambda,\mu,a}\Phi_p^{-1}=\mathbb{Q}_{\lambda,\mu,a}$.
Let us state some basic results on the homogeneous Poisson process, which will be useful in the sequel.
\subsection{Preliminary results}
Let $T=\left( T_k\right)_{k\geq 0}$ be a Poisson point process on $\mathbb{R}_{+}$ with intensity $\mu$. This means that $T_0=0$, and $\left( T_{k+1}-T_k,~k\ge 0\right)$ are i.i.d exponential r.v.'s with mean $1/\mu$. Let $\left(N_t,~t\ge 0 \right)$ be the counting process associated with $T$, that is $$\forall t\ge 0,~N_t=\sup \left\lbrace k\ge 0,~T_k\le t\right\rbrace.$$
The first result is well-known and elementary.
\begin{lem}\label{lemma21}
Let $M$ be a non negative random variable independent of $T$, and define \[R_M= \sup_{k\geq 0}\left\lbrace T_k; T_k\leq M\right\rbrace . \]
Then  $M-R_M\stackrel{(d)}=V\wedge M$ where $V$ and $M$ are independent, $V$ has an exponential distribution with mean $1/\mu$.\\
Moreover, on the event $\lbrace R_M>s\rbrace$, the conditional law of  $N_{R_M^-}-N_s$ given $R_M$ is Poisson with parameter $\mu(R_M-s)$.
\end{lem}
The second one is :
\begin{lem}\label{Poissoncombin}
Let $T=\left(T_k\right)_{k\geq 0}$ be a Poisson point process on
$\mathbb{R}_{+}$ with intensity $\mu$. $M$ a positive random
variable which is independent of $T$.
Consider the integer valued random variable $K$ such that $T_K=R_M$
and a second Poisson point process  $T^\prime=\left( T^\prime_k\right)_{k\geq
0}$ with intensity $\mu$, which is jointly independent of the first
and of $M$. Then $\overline{T}=\left( \overline{T}_k\right)_{k\geq 0}$ defined by:
\begin{displaymath}
\overline{T}_k=\left\{ \begin{array}{ll}
T_k & \textrm{if $k < K$}\\
T_K+T^\prime_{k-K+1} & \textrm{if $k\geq K$}
\end{array} \right.
\end{displaymath}
is a Poisson point process on $\mathbb{R}_+$ with intensity  $\mu$,
which is   independent of $R_M$.
\end{lem}
\bpf
Let  $\left( N_t, t\geq 0\right) $, $\left(\overline{N}_t,t\geq 0 \right)$ and
$\left({N}^\prime_t,t\geq 0 \right)$) be  the  counting processes   associated   to  $T$, $\overline{T} $ and ${T^\prime} $.\\
It suffices to prove that for  any  $n\ge  1,  \;  0<t_1<\dots<t_n$   and $k_1,\dots,k_n  \in\mathbb{N}^*$,
\[\xi_t=\mathbb{P}\Big(\overline{N}_{t_1}=k_1,\dots ,\overline{N}_{t_n}=k_n\vert R_M\Big)=e^{-\mu t_n}\prod_{i=1}^n\frac{(\mu (t_i-t_{i-1}))^{k_i-k_{i-1}}}{(k_i-k_{i-1})!}.\]
Since there is no harm in adding $t_i's$, we only need to do that computation on the event that there exists $2\leq i\leq n$ such that $t_{i-1}<R_M<t_i$, in which case a standard argument yields easily the claimed result, thanks to Lemma \ref{lemma21}. Indeed we have that
\begin{align*}
\xi_t &=\mathbb{P}\left(N_{t_1}=k_1,\cdots , N_{t_{i-1}}=k_{i-1},N_{R_M^-}+N^\prime_{t_i-t}=k_{i},\cdots ,N_{R_M^-}+N'_{t_n-t}=k_n\right) \\\\
&=\mathbb{P}\Bigg(N_{t_1}=k_1,\cdots , N_{t_{i-1}}-N_{t_{i-2}}=k_{i-1}-k_{i-2}, N_{R_M^-}-N_{t_{i-1}}
+N'_{t_i-R_M}=k_i-k_{i-1},\\ \\
 &\quad  N'_{t_{i+1}-R_M}-N'_{t_i-R_M}=k_{i+1}-k_i,\cdots , N'_{t_n-R_M}-N'_{t_{n-1}-R_M}=k_n-k_{n-1}\Bigg)\\\\
&= e^{-\mu t_n}\prod_{i=1}^n\frac{(\mu (t_i-t_{i-1}))^{k_i-k_{i-1}}}{(k_i-k_{i-1})!}.
\end{align*}
\epf
\subsection{Basic theorem}
We are now in a position to prove the next theorem, which says that the  tree  associated to the exploration process   $H^a$ defined in section \ref{prelim}  is a continuous-time binary Galton--Watson tree with death rate $\lambda$  and    birth  rate  $\mu$, killed at time $a$, and vice versa.
\begin{thm}\label{bijection}
For any $\lambda,\mu>0$ and $a\in (0,+\infty)$ $\left[\right.$including the case $a=+\infty$ in the case $\mu\leq\lambda\left. \right]$ ,
\[ \mathbb{Q}_{\lambda,\mu,a}=\mathbb{P}_{\lambda,\mu,a}\Phi_p^{-1}. \]
\end{thm}
\bpf The individuals making up the population represented by the tree whose law is $\mathbb{Q}_{\lambda,\mu,a}$, will be numbered: $\ell=1,2,...$. 1 is the ancestor of the whole family. The subsequent individuals will be identified below. We will show that this tree is explored by a process whose law is precisely $\mathbb{P}_{\lambda,\mu,a}$.
We introduce the  family   $(T_k^{\ell},k\geq 0, \ell \geq 1)$ of mutually independent Poisson processes with intensity $\mu$.  For any  $\ell \ge  1$, the process  $T_k^{\ell}$ describes the times of birth of the offsprings of the individual $\ell$. We define $U_\ell$ to be the lifetime of individual $\ell$. 
\begin{itemize}
 \item  $\textbf{Step 1}$: We start from the initial time $ t=0 $ and  climb up to the level $ M_1 $ of height $ U_1\wedge a $, where $ U_1 $ follows an exponential law with mean $ 1/\lambda$. We go down from $ M_1 $ until we find the most recent point of the Poisson process $(T_k^1) $ which gives the times of  birth of the offsprings of individual $ 1 $. So from Lemma 2.1, we have descended by $ V_1\wedge M_1 $, where $ V_1 $ follows an exponential law with mean $1/\mu$, and is independent of $M_1$. We hence reach the level $m_1=M_1-V_1\wedge M_1$.\\
If $m_1=0$, we stop, else we turn to
 \item $\textbf{Step 2}$: We give the label 2 to this last offspring of the individual 1, born at the time $m_1$. Let us define $(\bar{T}_k^2)$ by:
\begin{displaymath}
\bar{T}_k^2=\left\{ \begin{array}{ll}
T_k^1 & \textrm{if $k < K_1$}\\
T_{K_1}^1+T_{k-K_1+1}^2 & \textrm{otherwise}
\end{array} \right.
\end{displaymath}
where $K_1$ is such that $T_{K_1}^1=m_1$.\\
Thanks to Lemma \ref{Poissoncombin}, $(\bar{T}_k^2)$ is a Poisson process with intensity $\mu$ on $\mathbb{R}_+$, which is independent of $m_1$ and in fact also of $ (U_1, V_1)$.\\
 Starting from $m_1$, the exploration process climbs up to level $M_2=(m_1+U_2)\wedge a$, where $U_2$ is an exponential r.v. with mean $1/\lambda$, independent of $ (U_1, V_1) $. Starting from level $ M_2 $, we go down a height $ M_2\wedge V_2 $ where $V_2$ follows an exponential law with mean $1/\mu$ and is independent of $ (U_2, U_1, V_1) $, to find the most recent point of the Poisson process $(\bar{T}_k ^ 2) $. At this moment we are at the level $ m_2 = M_2-V_2\wedge M_2$. If $ m_2 = 0 $ we stop. Otherwise we give the label $3$ to the individual born at time $m_2$, and repeat step $2$ until we reach 0. See
 {\sc Figure 2}.
\end{itemize}
\begin{figure}[htbp] 
\centering
 \includegraphics[scale=0.6]{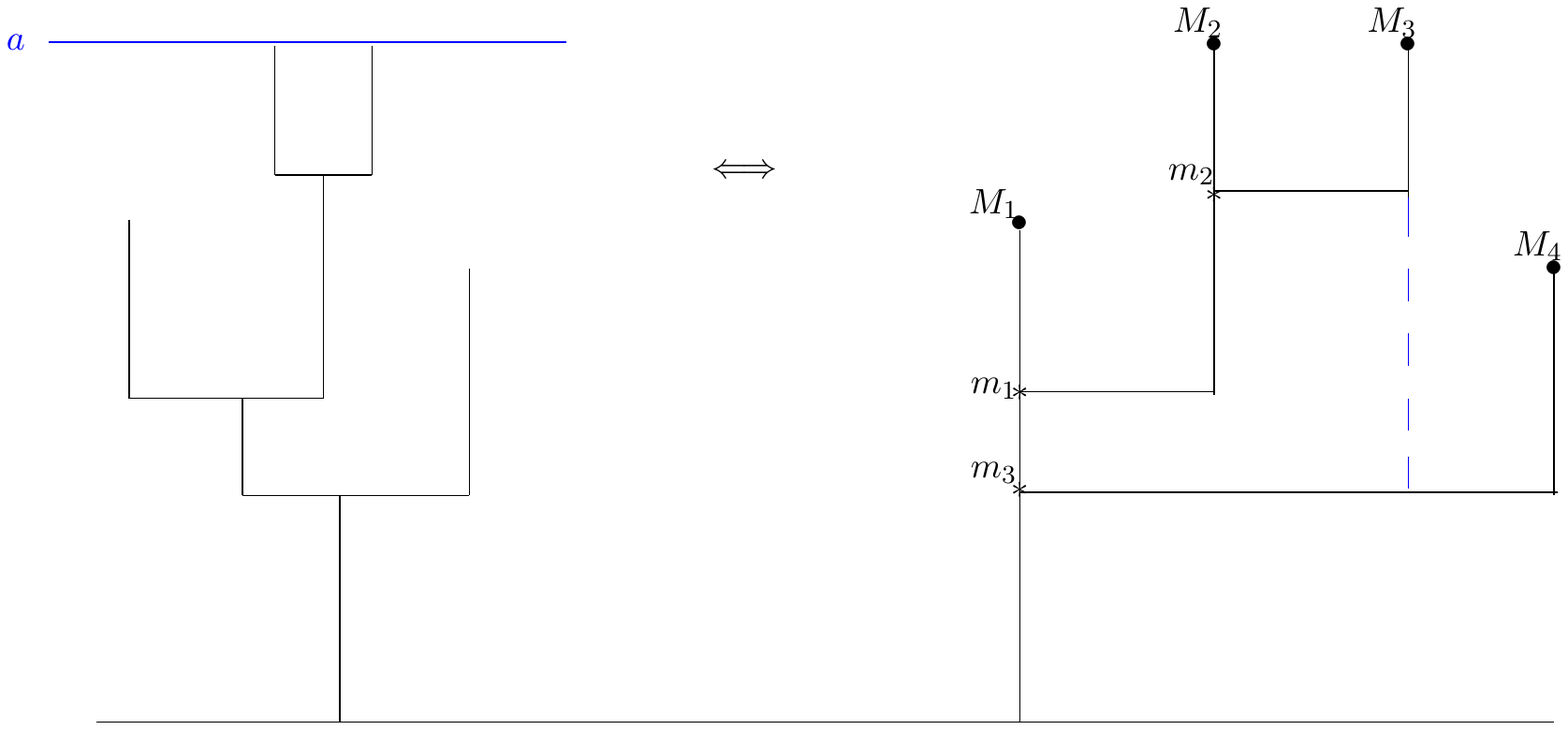}
  \caption {\label{Fig 2} Two equivalent ways of representing a binary tree}
\end{figure}    
Since either we have a reflection at level $a$ or $\mu\leq\lambda$, zero is reached $a.s.$ after a finite number of iterations. It is clear that the random variables $M_i $ $ $ and $ m_i $ determine fully the law $\mathbb{Q}_{\lambda,\mu,a}$ of the binary tree killed at time $t=a$ and they both have the same joint distribution as the levels of the successive local minima and maxima of the process $H^a$ under $\mathbb{P}_{\lambda,\mu,a}$. \epf
\subsection{A discrete Ray--Knight theorem}
For any $a,\mu,\lambda > 0$, we consider the exploration process $\{H^a_t,\ t\ge0\}$ defined in section \ref{prelim} which is reflected in the interval $[0,a]$ and stopped at the first moment it reaches zero for the $m$-th time. To this process we can associate a forest  of $m$ binary trees with birth rate $\mu$ and death rate $\lambda$, killed at time $t=a$, which all start with a single individual at the initial time $t=0$. Consider the branching process in continuous time $\left( Z^{a,m}_t,~t\geq 0\right)$ describing the number of offsprings alive at time $t$ of the $m$ ancestors born at time $0$, whose progeny is killed at time $t=a$. Every individual in this population, independently of the others, lives for an exponential time with parameter $\lambda$ and gives birth to offsprings according to a Poisson process of intensity $\mu$.
We now choose the slopes of the piecewise linear process $H^a$ to be $\pm2$ (i. e. $p=2$). 
We define the local time accumulated by $H^a$ at level $t$ up to time $s$:
\begin{equation}\label{loctime}
L^a_s(t)=\lim_{\eps\downarrow0}\frac{1}{\eps}\int_0^s{\bf1}_{\{t\le H^a_r<t+\eps\}}dr.
\end{equation}
$L^a_s(t)$ equals the number of pairs of branches of $H^a$ which cross the level $t$ between times $0$ and $s$. Note that a local minimum at level $t$ counts for two crossings, while a local maximum at level $t$ counts for none. We have the ``occupation time formula'': for any integrable function  $g$,
\[ \int_{0}^s g(H^a_r)dr=\int_0^{\infty} g(r)L^a_s(r)dr. \]
Let 
\begin{equation} \label{taua}
\tau^a_m=\inf \left\lbrace s >0: L^a_s(0)\ge m \right\rbrace .
\end{equation}
$L^a_{\tau_m}(t)$ counts the number of descendants of $m$ ancestors at time 0, which are alive at time $t$. Then we have
\begin{lem}\label{lemma2.6}
For all $\lambda,\mu>0$ and $a\in (0,+\infty)$ $\left[ \right.$including the case $a=+\infty$ in the case $\mu\leq\lambda\left.\right]$.
\[
\left\{L^a_{\tau^a_m}(t),\ t\geq 0,m\geq 1\right\}\equiv\left\{Z_t^{a,m},\
t\geq 0, m\geq 1\right\} ~~a.s..
\]
\end{lem}
We want now to establish a similar statement without the arbitrary parameter $a$. There remains a difficulty only in the supercritical case. For any $0<a<b$, we define the application $\Pi^{a,b}$  which maps continuous trajectories  with values in $[0,b]$ into trajectories with values in $[0,a]$ as follows. If $u\in C(\mathbb{R}_+,[0,b])$, 
$$\rho_u(s)=\int_0^s{\bf1}_{\lbrace u(s)\le a\rbrace}dr;~~\Pi^{a,b}(u)(s)=u(\rho^{-1}_u(s)).$$ 

\begin{lem}\label{lemma32}
\[ \Pi^{a,b}(H^b)\stackrel{(d)}=H^a \]
\end{lem}
\bpf
It is in fact sufficient to show that the conditional law of the level of the first local minimum of $H^b$ after crossing the level $a$ downwards, given  the past of $H^b$, is the same as the conditional law of the level of the first local minimum of $H^a$ after a reflexion at level $a$, given the past of $H^a$. This identity follows readily from the ``lack of memory" of the exponential law.  
\epf
This last Lemma says that reflecting under $a$, or chopping out the pieces of trajectory above level $a$, yields the same result (at least in law).

We now consider the case $p=2$. To each $\lambda,\mu >0$, $m\geq 1$, we associate the process $\left\lbrace Z^m_t,~t\geq 0\right\rbrace $ which describes the evolution of the number of descendants of $m$ ancestors, with birth rate $\mu$ and death rate $\lambda$. For each $a>0$ $\left[\right.$including the case $a=+\infty$ in the case $\mu\leq\lambda\left. \right]$, we let $\left(H^a_s,~s\geq 0 \right)$ denote the exploration process of the genealogical tree of this population killed at time $a$, $L^a$ denotes its local time and $\tau_m^a$ is defined by  \eqref{taua}. It follows readily from Lemma \ref{lemma32} that for any $0<a<b$,
\begin{align} \label{etoi}
\left( L^b_{\tau_m^b}(t),0\leq t<a, m\ge1\right)\stackrel{(d)}= \left( L^a_{\tau_m^a}(t),0\leq t<a, m\ge1\right).
\end{align}

The compatibility relation \eqref{etoi} implies the existence of a projective limit $\{\L_m(t),t\ge 0, m\ge1\}$ with values in $\mathbb{R}_+$, which is such that for each $a>0$,
\begin{equation}
\{\L_m(t),0 \leq t<a, m\ge1\}\stackrel{(d)}=\{L^a_{\tau^a_m}(t),0 \leq t<a, m\ge1\}. \label{comp}\end{equation}
We have the following ``discrete Ray--Knight Theorem''.
\begin{prop}\label{props}
$$\{\L_m(t), t\ge 0, m\geq 1\}\stackrel{(d)}=\left\{Z^m_t,\ t\ge0, m\geq 1\right\}.$$
\end{prop}
\bpf
It suffices to show that for any $a\ge0$, 
$$\{\L_m(t), 0\le t< a, m\ge1\}\stackrel{(d)}=\{Z^m_t, 0 \le t <a, m\ge1\}.$$
This result follows from \eqref{comp} and Lemma \ref{lemma2.6}.
\epf
\subsection{Renormalisation}
Let $x>0$ be arbitrary, and $N\geq 1$ be an integer which will eventually go to infinity. Let $\left( Z_t^{[Nx]}\right)_{t\geq 0}$ denote the branching  process which describes the number of descendants at time $t$ of $[Nx]$ ancestors, in the population with birth rate $\mu_N=\frac{\sigma^2}{2}N+\alpha$ and death rate $\lambda_N=\frac{\sigma^2}{2}N+\beta$, where $\alpha,\beta \geq 0$. We set  $$X^{N,x}_t =\frac{Z^{[Nx]}_t}{N}.$$
In particular we have that $X^{N,x}_0=\frac{[Nx]}{N}\longrightarrow x$ when $N\longrightarrow +\infty$. Let $H^{a,N}$ be the exploration process associated to $\left\lbrace  Z^{[Nx]}_t,~0\le t <a\right\rbrace $ defined in the same way as previously, but  with slopes $\pm 2N$, and  $\lambda$, $\mu$ are replaced by $\lambda_N$ and $\mu_N$. 
 We define also  $L_s^{a,N}(t)$, the local time accumulated by $H^{a,N}$ at level $t$ up to time $s$, as 
 \begin{align}\label{defLocTime} 
L^{a,N}_s(t)=\frac{4}{\sigma^2}\lim_{\eps\downarrow0}\frac{1}{\eps}\int_0^s{\bf1}_{\{t\le H^{a,N}_r<t+\eps\}}dr
\end{align}
The motivation of the factor $\frac{4}{\sigma^2}$ will be clear after we have taken the limit as $N\rightarrow\infty$.
$L^{a,N}_s(t)$ equals $4/N\sigma^2$ times the number of pairs of $t$-crossings of $H^{a,N}$  between times $0$ and $s$.
  Let 
  \begin{equation}\label{deftau}
   \tau_x^{a,N}=\inf \left\lbrace s >0: L_s^{a,N}(0)\ge \frac{4}{\sigma^2}\frac{[Nx]}{N} \right\rbrace.
  \end{equation}
We define for all  $N\ge 1$ the projective limit $\lbrace \L^{N}_x(t), t\ge0, x>0\rbrace$, which is such that for each $a>0$, 
$$\lbrace \L^{N}_x(t),0\le t< a, x>0\rbrace\stackrel{(d)}{=}\lbrace L^{a,N}_{\tau^{a,N}_x}(t),0\le t< a, x>0\rbrace.$$
Proposition \ref{props} translates as
\begin{lem}\label{RNd} We have the identity in law
$$\{\L^{N}_x(t), t\ge0, x>0\}\stackrel{(d)}{=}\left\{\frac{4}{\sigma^2}X^{N,x}_t, t\ge0, x>0\right\}.$$
\end{lem}

\section{Weak convergence}
\subsection{Weak convergence of $X^{N,x}$}
The following result describes the limit of the sequence of process $\left\lbrace X^{N,x},~N\geq 1\right\rbrace$, see e.g. Grimvall \cite{AG}.
\begin{prop}\label{convX}
$X^{N,x}\Rightarrow X^x$ as $N\to\infty$ for the topology of locally uniform convergence,
where $X^x$ is the unique solution of the following Feller SDE
$$X^x_t=x+(\alpha-\beta)\int_0^tX^x_rdr+\sigma\int_0^t\sqrt{X^x_r}dB_r,\ t\ge0.$$
\end{prop}
\subsection{Tightness criteria in $D([0,+\infty))$}
  Let us present a sufficient condition for tightness which will be useful below.
 Consider a sequence $\{X^n_t,\ t\ge0\}_{n\ge1}$ of one--dimensional semi--martingales, which is such that for each $n\ge1$,   
 \begin{align*}
  X^n_t&=X^n_0+\int_0^t\varphi^n_sds+M^n_t,\ 0\le t\le T;\\
  \langle M^n\rangle_t&=\int_0^t\psi^n_sds,\ t\ge 0;
  \end{align*}
  where for each $n\ge1$, $M^n_\cdot$ is a locally square--integrable martingale,  $\varphi^n$ and $\psi^n$ are progressively measurable processes with value in $\R$ and $\R_+$ respectively. Since our martingales $\lbrace M^n_t,~t\ge 0\rbrace$ will be discontinuous, we need to consider their trajectories as elements of $D\left(\left[0,+\infty\right)  \right) $, the space of right continuous functions with left limits at every point, from $\left[ 0,+\infty\right)$ into $\mathbb{R}$, which we equip with the Skorohood topology, see  Billingsley \cite{BI}.
  The following statement can be deduced from Theorems 16.10 and 13.4  in \cite{BI}.

  \begin{prop}\label{tightnessCriteria}
 A sufficient condition for the sequence $\{X^n_t,\ t\ge 0\}_{n\ge1}$
  to be tight in $D([0,\infty))$ is that both
  \begin{equation}\label{tightness0}
  \text{the sequence of r.v.'s }\{X^n_0,\ n\ge1\}\  \text{ is tight};
  \end{equation}
  and for  some  $c>0$,
  \begin{equation}\label{tightnessBound}
  \sup_{n\geq 1, s>0}\left( \vert \varphi^n_s \vert +\psi^n_s\right) \leq c.
  \end{equation}
  If moreover, for any $T>0$, as $n\to\infty$,
  $$\sup_{0\le t\le T}|M^n_t-M^n_{t^-}|\to0\quad\text{in probability},$$
  then any limit $X$ of a weakly converging subsequence of the original sequence $\{X^n\}_{n\ge1}$ is a. s. continuous.
  \end{prop}
\subsection{Tightness of $H^{a,N}$ \label{sec-tightHa}}
 Consider now the exploration process $\{H^{a,N}_s,\ s\ge0\}$ of the forest of trees representing the population $\{Z^{[Nx]}_t,\ 0\leq t<a \}$. Let $\{V^{a,N}_s,\ s\ge0\}$ be the $\{-1,1\}$--valued process which is such that $s$-a.e.
$
\frac{dH^{a,N}_s}{ds}=2NV^{a,N}_s
$.
The $\R_+\times\{-1,1\}$--valued process $\{(H^{a,N}_s,V^{a,N}_s),\ s\geq 0\}$ is a Markov process, which solves the following SDE~:
\begin{equation}\label{deriveH}
\begin{split}
\frac{dH_s^{a,N}}{ds}&=2NV^{a,N}_s,\quad H^{a,N}_0=0,  V^{a,N}_0=1,\\
dV^{a,N}_s&=2\mathbf{1}_{\lbrace V^{a,N}_{s^-}=-1\rbrace}dP^{+}_s-2\mathbf{1}_{\lbrace V^{a,N}_{s^-}=1\rbrace}dP^{-}_s+\frac{N
\sigma^2}{2} dL^{a,N}_s(0)- \frac{N \sigma^2}{2}dL^{a,N}_s(a^-),
\end{split}
\end{equation}
where $\{P^+_s,\ s\ge0\}$ and $\{P^-_s,\ s\ge0\}$ are two mutually
independent Poisson processes, with intensities respectively $$\sigma^2 N^2+2\alpha N\quad\text{ and } \quad
\sigma^2N^2+2\beta N,$$
$L^{a,N}_s(0)$ and $L^{a,N}_s(a^-)$ denote respectively the number of visits to $0$ and $a$ by the process $H^{a,N}$ up to time $s$, multiplied by $4/N\sigma^2$ (see \eqref{defLocTime}). These two terms in the expression of $V^{a,N}$ stand for the reflection of $H^{a,N}$ above $0$ and below $a$. Note that our definition of $L^{a,N}$ makes the mapping $t\longrightarrow L^{a,N}_s(t)$ right continuous for each $s>0$. Hence $L^{a,N}_s(t)=0$ for $t\geq a$, while $L^{a,N}_s(a^-)=\lim_{n\rightarrow \infty}L^{a,N}_s(a-\frac{1}{n})>0$ if $H^{a,N}$ has reached the level $a$ by time $s$.

We now write a sub--martingale problem satisfied by the process  $\{(H^{a,N}_s,V^{a,N}_s),$ $\ s\geq 0\}$. We are not interested in writing it for arbitrary functions of the two variables $(h,v)$, but rather for specific functions, which will be convenient for taking the limit as $N\to\infty$. Note that the process $\{V^{a,N}_s,\ s\ge0\}$ oscillates faster and faster as $N$ grows, and that in the limit some averaging takes place.
We thus implement the so called ``perturbed test function method'' used in stochastic averaging, see e. g. Ethier and Kurtz \cite{EK}.
For $f\in C^2(\R)$, let
\begin{align*}
f^N(h,v)&=f(h)+\frac{v}{N\sigma^2}f'(h),\\
A^Nf^N(h,v)&=\frac{2}{\sigma^2}f''(h)+{\bf1}_{\{v=-1\}}\frac{4\alpha}{\sigma^2}f'(h)
-{\bf1}_{\{v=+1\}}\frac{4\beta}{\sigma^2} f'(h).
\end{align*}
If $f'(0)\ge 0$ and $f'(a)\leq 0$, then 
\begin{equation}\label{subMa}
M^{f,N,a}_s:=f^N(H^{a,N}_s,V^N_s)-f^N(0,1)-\int_0^s A^Nf^N(H^N_r,V^N_r)dr
\end{equation}
is a local sub--martingale. If we rather choose successively $f(h)=h$ and $f(h)=h^2$, we deduce from \eqref{deriveH} that there exist two local martingales  $\{M^{1,a,N}_s,\ s\ge0\}$ and $\{M^{2,a,N}_s,\ s\ge0\}$ such that
\begin{equation}\label{eqapprox}
\begin{split}
H^{a,N}_s+\frac{V^{a,N}_s}{N\sigma^2}&=\frac{1}{N\sigma^2}+
\frac{4\alpha}{\sigma^2}\int_0^s{\bf1}_{\{V^{a,N}_r=-1\}}dr-
\frac{4\beta}{\sigma^2}\int_0^s{\bf1}_{\{V^{a,N}_r=+1\}}dr\\
&\quad+\frac{1}{2}\left[L^{a,N}_s(0)-L^{a,N}_{0^+}(0)\right]-\frac{1}{2}L^{a,N}_s(a^-)+M^{1,a,N}_s,
\end{split}
\end{equation}
\begin{align*}
(H^{a,N}_s)^2+\frac{2}{N\sigma^2}H^{a,N}_s V^{a,N}_s &=\frac{4}{\sigma^2}s+
\frac{8\alpha}{\sigma^2}\int_0^s{\bf1}_{\{V^{a,N}_r=-1\}}H^{a,N}_rdr\\
&\quad-\frac{8\beta}{\sigma^2}\int_0^s{\bf1}_{\{V^{a,N}_r=+1\}}H^{a,N}_rdr-aL^{a,N}_s(a^-)+M^{2,a,N}_s.
\end{align*}
It follows from the above computations that
\begin{align} \label{crochHa}
\langle M^{1,a,N}\rangle_s &=\frac{4}{\sigma^2}s+\frac{8\alpha}{N\sigma^4}
\int_0^s {\bf1}_{\{V^{a,N}_{r}=-1\}}dr+\frac{8\beta}{N\sigma^4}
\int_0^s {\bf1}_{\{V^{a,N}_{r}=1\}}dr ,
\end{align}
and from \eqref{crochHa} that $\{M^{1,a,N}_s,\ s\ge0\}$ is in fact a martingale.  One difficulty which we want to get rid of is the local time terms in the expression for $H^{a,N}_s+\frac{V^{a,N}_s}{N\sigma^2}$, which introduce some additional complication for checking tightness. For that sake, we consider a new pair of processes $(G^{a,N},W^{a,N})$, which is $\R\times\{-1,1\}$--valued and satisfies:
\begin{align*}
G^{a,N}_s&=2N\int_0^s W^{a,N}_r dr,\\
W^{a,N}_s&=1+\sum_{i\in\mathbb{Z}}\Bigg\lbrace 2\int_0^s{\bf1}_{\lbrace ai\leq G^{a,N}_r\leq (i+1)a\rbrace}(-1)^{i}{\bf1}_{\lbrace W^{a,N}_{r^-}=-(-1)^{i}\rbrace}dP^+_r\\
&\quad-2\int_0^s{\bf1}_{\lbrace ai\leq G^{a,N}_r\leq (i+1)a\rbrace}(-1)^{i}{\bf1}_{\lbrace W^{a,N}_{r^-}=-(-1)^{i}\rbrace}dP^-_r\Bigg\rbrace
\end{align*}
with the same $P^+$ and $P^-$ as above. We claim that $$H^{a,N}=a.s.~\lim_{k \rightarrow\infty}\varphi_k(G^{a,N}), $$ 
$$V^{a,N}=\sum_{i\in\mathbb{Z}}(-1)^{i} {\bf1}_{\lbrace ai\leq G^{a,N}\leq (i+1)a\rbrace}W^{a,N},$$
where 
\begin{align*}
\varphi_k=\psi_k\circ\cdots\circ \psi_1,
\end{align*}
and for every $j$,  the mapping $\psi_j$ from $\mathbb{R}$ into $\mathbb{R}$ is defined by:
\begin{displaymath}
\psi_j(x)=\left\{ \begin{array}{ll}
|x|, & \textrm{ if $j$ is odd;}\\
a-|x-a|, & \textrm{ if $j$ is even.}
\end{array} \right.
\end{displaymath}
Indeed, since $G^{a,N}_s$ is locally bounded, to each $r>0$ we can associate a random index $k$ such that $0\leq \psi_k(G^{a,N}_s)\leq a$, for $0\leq s\leq r$, which implies that $\phi_{k+j}(G^{a,N}_s)=\phi_k(G^{a,N}_s)$, for $0\leq s\leq r$, $j\geq 1$. 
Note that $\psi_1$ consists in reflecting the $G^{a,N}$ trajectory above 0, $\psi_2$ reflecting below $a$. Those operations are repeated until the thus obtained trajectory stays in $[0,a]$. The reader can convince himself that it then coincides with $H^{a,N}$.
Tightness of $\{G^{a,N}\}$ will imply that of  $\{H^{a,N}\}$, since 
\[ \forall s, t ~~~|H^{a,N}_s-H^{a,N}_t| \leq  |G^{a,N}_s-G^{a,N}_t|.\] 
We have 
\begin{align}
G^{a,N}_s+\frac{W^{a,N}_s}{N\sigma^2}&=\frac{1}{N\sigma^2}+
\frac{4\alpha}{\sigma^2}\sum_{i\in \mathbb{Z}}\int_0^s{\bf1}_{\lbrace ai\leq G^{a,N}_r\leq (i+1)a\rbrace}(-1)^{i}{\bf1}_{\lbrace W^{a,N}_{r^-}=-(-1)^{i}\rbrace}dr\nonumber
\\
&\qquad-
\frac{4\beta}{\sigma^2}\sum_{i\in \mathbb{Z}}\int_0^s{\bf1}_{\lbrace ai\leq G^{a,N}_r\leq (i+1)a\rbrace}(-1)^{i}{\bf1}_{\lbrace W^{a,N}_{r^-}=(-1)^{i}\rbrace}dr+\tilde{M}^{1,N,a}_s\label{Ga}
\\
\langle \tilde{M}^{1,N,a}\rangle_s &=\frac{4}{\sigma^2}s+\frac{8\alpha}{N\sigma^4}\sum_{i\in \mathbb{Z}}\int_0^s {\bf1}_{\lbrace ai\leq G^{a,N}_r\leq (i+1)a\rbrace}{\bf1}_{\lbrace W^{a,N}_{r^-}=-(-1)^{i}\rbrace}dr\nonumber \\ &\qquad+\frac{8\beta}{N\sigma^4}\sum_{i\in \mathbb{Z}}
\int_0^s {\bf1}_{\lbrace ai\leq G^{a,N}_r\leq (i+1)a\rbrace}{\bf1}_{\lbrace W^{a,N}_{r^-}=(-1)^{i}\rbrace}dr.\label{Ma}
\end{align}
From \eqref{Ga}, \eqref{Ma} and Proposition \ref{tightnessCriteria} follows tightness of the lefthand side of \eqref{Ga}. Since moreover $N^{-1}W^{a,N}_s\rightarrow 0$ a.s. uniformly with respect to $s$, the sequence $\{G^{a,N},~N\ge1\}$ is tight in $D([0,+\infty))$. Because $G^{a,N}$ is a. s. continuous for each $N\ge1$, it follows from a well known property of Skorohod's topology~:
\begin{lem}\label{tightHa}
For any $a>0$, the sequence $\{H^{a,N}_s,\ s\ge0\}_{N\ge1}$ is tight in $C([0,\infty))$.
\end{lem}
\begin{rem}
In the subcritical case ($\alpha\leq\beta$), we can choose $a=+\infty$, which simplifies the above construction. $H^N$ is obtained from $G^N$ by reflection above 0 ($H^N\equiv |G^N|$), and $G^N$ is defined by:
\begin{align*}
G^N_s&=2N\int_0^s W^N_r dr,\\
W^N_s&=1+2\int_0^s\sign(G^N_{r}){\bf1}_{\{W^N_{r^-}=-\sign(G^N_r)\}}dP^+_r-2\int_0^s\sign(G^N_{r}){\bf1}_{\{W^N_{r^-}=\sign(G^N_r)\}}dP^-_r
\end{align*}
\end{rem}
\subsection{Weak convergence of $H^{a,N}$}
Let us state our convergence result.
 \begin{thm}\label{convHa} For any $a>0$ $\left[\right.$including the case $a=+\infty$ in the case $\alpha\leq\beta\left. \right]$,
$H^{a,N}\Rightarrow H^a$ in $C([0,\infty))$ as $N\to\infty$, where $\{H^a_s,\ s\ge0\}$ is the process $$\frac{2(\alpha-\beta)}{\sigma^2}s+\frac{2}{\sigma}B_s$$ reflected in $[0,a]$. In other words, $H^a$ is the unique weak solution of the reflected SDE
 \begin{equation}\label{browref}
 H^a_s=\frac{2(\alpha-\beta)}{\sigma^2}s+\frac{2}{\sigma}B_s+\frac{1}{2}L_s(0)-\frac{1}{2}L_s(a^-).
 \end{equation}
 \end{thm}
 The statement that $\{H^a_s,\ s\ge0\}$ is the process $\left( \frac{2(\alpha-\beta)}{\sigma^2}s+\frac{2}{\sigma}B_s,~s\geq 0\right)$ reflected in $[0,a]$ amounts to saying (see Stroock and Varadhan \cite{SV}) that for any $f\in C^2(\mathbb{R})$ with $f'(0)\geq0$, $f'(a)\leq 0$,
\begin{equation*}
 M^{f}_s:=f(H^a_s)-f(H^a_0) -\frac{2}{\sigma^2}\int_0^s \left[  f''(H^a_r)-(\alpha-\beta)f'(H^a_r)\right]dr,
\end{equation*}
is a sub--martingale. It remains to establish this property by taking the weak limit in \eqref{subMa}. This will follow readily from
\begin{lem}\label{convVHa}
For any sequence $\left( U^N,N\ge 1\right)\subset C([0,+\infty))$ which is such that $U^N\Rightarrow U$ as $N\rightarrow\infty$, for all $s>0$,
$$\int_0^s  {\bf1}_{\{V^{a,N}_r=1\}}U^N_rdr\Rightarrow\frac{1}{2}\int_0^s U_rdr,\quad \int_0^s  {\bf1}_{\{V^{a,N}_r=-1\}}U^N_rdr\Rightarrow\frac{1}{2}\int_0^sU_rdr.$$
\end{lem}
\bpf
 It is an easy exercise to check that the mapping
 $$\Phi : C([0,+\infty))\times C_{\uparrow}([0,+\infty))\to C([0,+\infty))
 $$
 defined by
 $$\Phi(x,y)(t)=\int_0^tx(s)dy(s),$$
 where $C_{\uparrow}([0,+\infty))$ denotes the set of increasing continuous functions from $[0,\infty)$ into $\R$, and the three spaces are equipped with the topology of locally uniform convergence,  is continuous. Consequently it suffices to prove that
 locally uniformly in $s>0$,
 $$\int_0^s {\bf1}_{\{V^{a,N}_r=1\}}dr\to\frac{s}{2}$$ in probability, as $N\to\infty$. In fact since both the sequence of processes and the limit are continuous and monotone, it follows from an argument ``\`a la Dini'' that  it suffices to prove 
 \begin{lem}\label{convV}
 For any $s>0$,
 $$\int_0^s {\bf1}_{\{V^{a,N}_r=1\}}dr\to\frac{s}{2}$$ in probability, as $N\to\infty$.
 \end{lem} 
\bpf
We have (the second line follows from \eqref{deriveH})
\begin{align*}
\int_0^s {\bf1}_{\{V^{a,N}_r=1\}}dr+\int_0^s {\bf1}_{\{V^{a,N}_r=-1\}}dr&=s,\\
\int_0^s {\bf1}_{\{V^{a,N}_r=1\}}dr-\int_0^s {\bf1}_{\{V^{a,N}_r=-1\}}dr&=(2N)^{-1}H^{a,N}_s.
\end{align*}
From Lemma \ref{tightHa} follows readily that $(2N)^{-1}H^{a,N}_s\to0$ in probability, as $n\to\infty$.
We conclude by adding the two above identities.
\epf
\begin{cor}\label{corro}
For each $a>0$ $\left[\right.$including the case $a=+\infty$ in the case $\alpha\leq\beta\left. \right]$, 
$$\left(H^{a,N},M^{1,N,a},L^{a,N}_{\centerdot}(0),L^{a,N}_{\centerdot}(a^-) \right)\Longrightarrow \left( H^a,\frac{2}{\sigma}B, L^a_{\centerdot}(0), L^a_{\centerdot}(a^-)\right), $$
where $B$ is as above, $L^a_{\centerdot}(0)$ $\left( \right.$resp. $L^a_{\centerdot}(a^-)\left.\right)$ denotes the local time of the continuous semi--martingale $H^a$ at level $0$ $\left(\right.$resp. at level $a^-\left. \right) $. 
\epf
\end{cor}
\bpf Concerning tightness, we only need to take care of the third and fourth terms in the quadruple.
Consider the function $f^N(h,v)$ associated to some $f\in C^2(\R)$ such that $f'(0)=1$ and $f'(a)=0$. 
We deduce from \eqref{deriveH}
\begin{equation}\label{tempsloc}
\begin{split}
L^{a,N}_s(0)&=2f(H^{a,N}_s)+2\frac{V^{a,N}_s}{N\sigma^2}f'(H^{a,N}_s)-2f(0)-\frac{2}{N\sigma^2}f'(0)
-\frac{4}{\sigma^2}\int_0^sf''(H^{a,N}_r)dr\\&-\frac{8}{\sigma^2}\int_0^sf'(H^{a,N}_r)(\alpha{\bf1}_{\{V^N_r=-1\}}-\beta{\bf1}_{\{V^N_r=1\}})dr
-2M^{f,N}_s-2\tilde{M}^{f,N}_s,
\end{split}
\end{equation}
where $M^{f,N}$ and $\tilde{M}^{f,N}$ are martingales such that
$$
\langle M^{f,N}\rangle_s= \frac{4}{\sigma^2}\int_0^s[f'(H^{a,N}_r)]^2dr ,\quad \langle \tilde{M}^{f,N}\rangle_s\le\frac{c(f)}{N}s. $$
 Tightness of the local time terms follows from this formula (and a similar expression for $L^{a,N}_s(a^-)$).
 Then $\big(H^{a,N}, M^{1,N,a},L^{a,N}_{\centerdot}(0),$ $L^{a,N}_{\centerdot}(a^-)\big)_{N\ge 1}$ is tight in $C([0,\infty))\times\left[D([0,\infty))\right]^3$.

Moreover any weak limit of $M^{1,N,a}$ along a subsequence equals $\frac{2}{\sigma}B$, since \\
$<M^{1,N,a}>_s\rightarrow \frac{4}{\sigma^2}s$ and the jumps of $M^{1,N,a}$ are equal in amplitude to $\frac{2}{N\sigma^2}$. It then follows by taking the limit in \eqref{tempsloc} (and in a similar formula for $L^{a,N}_s(a^-)$) that any weak limit of 
$\left(H^{a,N},M^{1,N,a},1/2 L_s(0),L^{a,N}_{\centerdot}(a^-) \right)$ along a converging subsequence takes the form $\left( H^a,\frac{2}{\sigma}B, 1/2 L_s(0),1/2 L_s(a^-)\right)$.  Finally from Theorem $5.7$ in \cite{SV} the limit is unique, hence the whole sequence converges.
\epf

\section{Generalized Ray Knight Theorem}
In this section we give an new proof of Delmas' generalization of the second Ray Knight Theorem. Define $L_\cdot(0)$ be the local time of $H$ at level $0$, and in the
subcritical case $\alpha\le\beta$
\[ \tau_x=\inf\{  s>0; L_s(0)>\frac{\sigma^2}{4}x\} . \]
In the supercritical case, of course the construction is more complex. It follows from Lemma \ref{lemma32} and Corollary \ref{corro} (see also Lemma 2.1 in \cite{JFD}) that for any $0<a<b$, 
\begin{equation}\label{lemmeJFD}
\Pi^{a,b}(H^b)\stackrel{(d)}=H^a,
\end{equation}
where $H^a$ [resp. $H^b$] is Brownian motion multiplied by $2/\sigma$, with
drift $2(\alpha-\beta)s/\sigma^2$, reflected in the interval $[0,a]$ [resp. $[0,b]$], see Theorem \ref{convHa}. 
Now define for each $a,x>0$,
$$\tau^a_x=\inf\{s>0,\ L^a_s(0)>\frac{4}{\sigma^2}x\}.$$
It follows from \eqref{lemmeJFD} that, as in the discrete case, $\forall ~0<a<b$,
$$\{L^b_{\tau_x^b}(t),\ 0\le t< a, x>0\}\stackrel{(d)}=\{L^a_{\tau_x^a}(t),\ 0\le t<a, x>0\}.$$
Consequently we can define the projective limit, which is a process $\{\L_x(t),\ t\ge0, x>0\}$ such that for each $a>0$,
$$\{\L_x(t),\ 0\le t<a, x>0\}\stackrel{(d)}=\{L^a_{\tau_x^a}(t),\ 0\le t<a, x>0\}.$$
 We have the (see Theorem 3.1 in Delmas \cite{JFD} )
\begin{thm}[Generalized Ray Knight theorem]
\[ \{\L_x(t), t\ge0,x>0\} \stackrel{(d)}=\left\{\frac{4}{\sigma^2}X^x_t,t\ge 0,x>0 \right\},\]
where $X^x$ is the Feller branching diffusion process, solution of the SDE
$$X^x_t=x+(\alpha-\beta)\int_0^tX^x_rdr+\sigma\int_0^t\sqrt{X^x_r}dB_r,\ t\ge0.$$
\end{thm}
\bpf
Since both sides have stationary independent increments in $x$, it suffices to show that for any $x>0$,
\[ \{\L_{x}(t),\ t\ge0\}\stackrel{(d)}=\left\{\frac{4}{\sigma^2}X^x_t,\ t\ge0\right\}.\]
Fix an arbitrary $a>0$.
By applying the ``occupation time formula'' to $H^{a,N}$, and Lemma \ref{RNd}, we have for any
$g\in C(\R_+)$ with support in $[0,a]$,
\begin{align}\label{a}
\frac{4}{\sigma^2}\int_0^{\tau^{a,N}_x} g(H^{a,N}_r)dr &=
\int_0^\infty g(t)L^{a,N}_{\tau_x^{a,N}}(t)dt \nonumber \\
&\stackrel{(d)}=\frac{4}{\sigma^2}\int_0^\infty g(t)X^{N,x}_tdt
\end{align}
We deduce clearly from Proposition \ref{convX}
\begin{equation}\label{b}
\int_0^\infty g(t)X^{N,x}_t dt\Longrightarrow\int_0^\infty g(t)X^x_tdt.
\end{equation}
Let us admit for a moment that as $N\rightarrow \infty$
\begin{equation} \label{c}
\int_0^{\tau^{a,N}_x} g(H^{a,N}_r)dr\Longrightarrow\int_0^{\tau^a_x} g(H^a_r)dr
\end{equation}

From the occupation time formula for the continuous semi-martingale $(H^a_s,s\ge 0)$, we have that
\begin{equation}\label{d}
\frac{4}{\sigma^2}\int_0^{\tau^a_x}g(H^a_r)dr=\int_0^{\infty}g(t)L^a_{\tau^a_x}(t)dt. 
\end{equation}

We deduce from \eqref{a}, \eqref{b}, \eqref{c} and \eqref{d} that for any
 $g\in C(\mathbb{R}_+)$ with compact support in $[0,a]$, 
\[\frac{4}{\sigma^2}\int_0^{\infty} g(t)X^x_tdt\stackrel{(d)}=\int_0^\infty g(t)\L_{x}(t)dt.\]

In fact, this same argument can be slightly generalized, proving that for any $n$, any $g_1,\cdots,g_n$ with compact support, we have the following equality in distribution
\begin{equation*}
\left(\frac{4}{\sigma^2}\int_0^{\infty} g_1(t)X^x_tdt,\cdots,\frac{4}{\sigma^2}\int_0^{\infty} g_n(t)X^x_tdt \right)\stackrel{(d)}=\left( \int_0^{\infty} g_1(t)\L_{x}(t)dt,\cdots,\int_0^{\infty} g_n(t)\L_{x}(t)dt\right).  
\end{equation*}
Since both processes $\left( X^x_t,t\geq 0\right) $ and $\left( \L_{x}(t),t\geq 0\right) $ are $a.s.$ continuous, the theorem is proved.
 \epf

It remains to prove \eqref{c}, which clearly is a consequence of (recall the definition \eqref{deftau} of $\tau^N_x$)
\begin{prop}
\[ (H^{a,N},\tau^{a,N}_x)\Longrightarrow (H^a,\tau^a_x). \]
\end{prop}
\bpf For the sake of simplifying the notations, we suppress the superscript $a$.
Let us define the  function $\phi$ from $\mathbb{R}_+\times C_{\uparrow}([0,+\infty))$ into $\mathbb{R}_+$  by
\[ \phi(x,y)=\inf \lbrace s>0: y(s)>\frac{4}{\sigma^2}x \rbrace .\]
For any fixed $x$, the function $\phi(x,.)$ is continuous in the neighborhood of a function $y$ which is strictly increasing at the time when it first reaches the value $4x/\sigma^2$.
 Define
 \[\tau'^N_x:=\phi(x,L^N_.(0)).\]
 We note that for any $x>0$, $s\longmapsto L_s(0)$ is a.s. strictly increasing at time $\tau_x$, which is a stopping time. This follows from the strong Markov property, the fact that $H_{\tau_x}=0$, and $L_\epsilon(0)>0$, for all $\epsilon>0$.
 Consequently $\tau_x$ is $a.s.$ a continuous function of the trajectory $L_.(0)$, and from Corollary \ref{corro}
  \[ (H^N,\tau'^N_x)\Longrightarrow (H,\tau_x). \]
 It remains to prove that $\tau'^N_x-\tau^N_s\longrightarrow 0$ in probability. For any $y<x$, for $N$ large enough
 \[ 0\le \tau'^N_x-\tau^N_x\le \tau'^N_x -\tau'^N_y.\]
 Clearly $ \tau'^N_x -\tau'^N_{y}\Longrightarrow  \tau_x -\tau_{y},$ hence for any $\varepsilon>0$,
 \[0\le \limsup_{N}\mathbb{P}\big(\tau'^N_x -\tau^N_x\ge\varepsilon\big)\le \limsup_N\mathbb{P}\big(\tau'^N_x -\tau'^N_y\ge\varepsilon\big)\le\mathbb{P}\big(\tau_x -\tau_y\ge\varepsilon\big).  \]
The result follows, since $\tau_y\rightarrow\tau_{x^-}$ as $y \rightarrow x$, $y<x$, and
$\tau_{x^-}=\tau_x$ a. s.
\epf

\noindent{\bf Acknowledgements} This work was partially supported by the ANR MANEGE (contract ANR-09-BLANC-0215-03).
The authors wish to thank an anonymous Referee, whose excellent and detailed reports helped us to correct one error and improve significantly the presentation of our results.

 \end{document}